\documentclass[graybox]{svmult}

\usepackage{mathptmx}       
\usepackage{helvet}         
\usepackage{courier}        
\usepackage{type1cm}        
\usepackage{makeidx}         
\usepackage{graphicx}        
\usepackage{multicol}        
\usepackage[bottom]{footmisc}
\usepackage{amssymb,xcolor}

\newcommand{\be}{\begin{equation}}
\newcommand{\ba}{\begin{array}}
\newcommand{\ee}{\end{equation}}
\newcommand{\ea}{\end{array}}
\newcommand{\Exp}{{\rm I\hspace{-0.8mm}E}}

\def\N{\mathbb{N}}
\def\Z{\mathbb{Z}}
\def\R{\mathbb{R}}
\newcommand{\eqref}[1]{(\ref{#1})}%
\newtheorem{assumption}{Assumption}{\bfseries}{\rmfamily}
\makeindex  

\begin{document}

\title*{Zero-range process in random environment}
\author{C. Bahadoran, T.S. Mountford, K. Ravishankar and E. Saada}
\institute{
C. Bahadoran \at Laboratoire de Math\'ematiques Blaise Pascal, 
Universit\'e Clermont Auvergne, 63177 Aubi\`ere, France. 
\email{Christophe.Bahadoran@uca.fr}
\and T.S. Mountford \at Institut de Math\'ematiques, 
\'Ecole Polytechnique F\'ed\'erale, Lausanne, Switzerland.
\email{thomas.mountford@epfl.ch}
\and K. Ravishankar \at NYU-ECNU Institute of Mathematical Sciences at NYU Shanghai,
  3663 Zhongshan Road North, Shanghai, 200062, China.
\email{ ravi101048@gmail.com }
\and E. Saada \at CNRS, UMR 8145, MAP5,
Universit\'e de Paris, Campus Saint-Germain-des-Pr\'es, France.
\email{Ellen.Saada@mi.parisdescartes.fr}
}
\maketitle

\abstract*{
We survey our recent articles dealing with one dimensional attractive 
zero range processes moving
under site disorder.  We suppose that the underlying random walks 
are biased to the right and so hyperbolic scaling is expected.  
Under the conditions of our model the process admits a 
 maximal invariant measure.   The initial focus of the project 
 was to find conditions on
 the initial law to entail convergence in distribution to this maximal distribution,
 when it has a finite density.  
 Somewhat surprisingly, necessary and sufficient conditions were found. 
 In this part hydrodynamic results were employed chiefly as a tool to show distributional convergence but 
 subsequently  we developed  a theory for hydrodynamic limits  treating 
 profiles possessing densities that did not admit corresponding equilibria.
  Finally we derived strong local equilibrium results. 
}

\abstract{
We survey our recent articles dealing with one dimensional attractive 
zero range processes moving
under site disorder.  We suppose that the underlying random walks 
are biased to the right 
and so hyperbolic scaling is expected.  Under the conditions 
of our model the process admits a 
 maximal invariant measure.   The initial focus of the project 
 was to find conditions on
 the initial law to entail convergence in distribution to this maximal distribution,
 when it has a finite density.  
 Somewhat surprisingly, necessary and sufficient conditions were found. 
 In this part hydrodynamic results were employed chiefly as a tool 
 to show distributional convergence but 
 subsequently  we developed  a theory for hydrodynamic limits  treating 
 profiles possessing densities that did not admit corresponding equilibria.
  Finally we derived strong local equilibrium results. 
}

\section{Introduction}
\label{sec_intro}
The  asymmetric zero-range processs (AZRP) with site disorder was introduced in \cite{ev}
to study condensation phenomena. It is a conservative interacting particle system 
whose dynamics is determined by a  jump rate  function $g:\mathbb{N}\to\mathbb{N}$, a function 
 $\alpha :\mathbb{Z}^d\to\mathbb{R}_+$ (called the \textit{environment} or \textit{disorder}),
and a jump distribution $p(.)$ on $\mathbb{Z}^d$, for $d\ge 1$.    
A particle leaves site $x$ at rate $\alpha(x)g[\eta(x)]$, where $\eta(x)$ denotes the 
current number of particles at $x$, and moves to $x+z$, where  site  $z$ is chosen at random 
with distribution $p(.)$. As explained later on, this  model has a whole family of 
product invariant measures carrying different mean densities;
it exhibits a {\em critical density}  $\rho_c$, i.e.  no
product invariant measure exists above $\rho_c$ (\cite{fk,bfl}), 
if the function $g$ is bounded,
$\alpha$ has averaging properties plus a proper tail assumption. 
This can be interpreted as a {\em phase transition}. \par\medskip 
 In this review paper, we consider the one-dimensional attractive 
 nearest neighbour process, that is
  $d=1$, $p(1)+p(-1)=1$ and $g$ nondecreasing. 
  We summarize the papers \cite{bmrs1,bmrs2,bmrs3,bmrs4},
  by giving  their results and the main ideas of their proofs. 
In these papers,  we  developed robust approaches to study 
various aspects of the phase transition mentioned above.\par
One aspect is the {\em mass escape} phenomenon.
Suppose the process is started from 
a given configuration where the global empirical 
density of particles is greater than $\rho_c$. 
One usually expects convergence to 
to the extremal invariant measure carrying the same density as the initial state. 
However, in this case such a measure does not exist. 
When  $g(n)=\min(n,1)$  and $p(1)=1$ it was shown in \cite{afgl} 
that the system converges to the maximal invariant measure (thereby implying 
a loss of mass). This was established 
in \cite{bmrs1,bmrs2} for the general nearest neighbour model under a weak 
convexity assumption, and we showed that
this could fail for non nearest neighbour jump kernels.\par
Phase transition also arises  in the hydrodynamic limit. We show in \cite{bmrs3} that
the hydrodynamic behaviour of our process is given under hyperbolic time scaling 
by entropy solutions of a scalar conservation law
\be\label{burgers_intro}
\partial_t\rho(t,x)+\partial_x[f(\rho(t,x))]=0
\ee
where $\rho(t,x)$ is the local particle density field, with a  macroscopic flux  
function $\rho\mapsto f(\rho)$ that is increasing up to $\rho_c$ and constant thereafter.
\par
The natural question following hydrodynamic limit is that of {\em local equilibrium}. 
In general, for a conservative particle system endowed with a family 
$(\nu_\rho)_{\rho
}$ of extremal invariant measures
denotes the set of allowed macroscopic densities), 
the local equilibrium property 
states that the distribution of the microscopic particle configurations 
 at a macroscopic time $t\geq 0$  around a site 
with macroscopic location $x\in\mathbb{R}$ is close 
to $\nu_{\rho(t,x)}$, where $\rho(t,x)$ 
is the hydrodynamic density, here given by (\ref{burgers_intro}). 
This property has a 
weak (space-averaged) and a strong (pointwise) formulation, see e.g. \cite{kl}.
But the local equilibrium property is  expected to be 
 {\em wrong} at \textit{supercritical} hydrodynamic densities,  that is, such that  
 $\rho(t,x)>\rho_c$, since a corresponding 
 equilibrum measure does not exist. 
 This already poses a problem at the level of the hydrodynamic 
 limit (\cite{bmrs3}),
since in the usual heuristic for (\ref{burgers_intro}), 
the macroscopic flux function $f$ is the expectation of  
microscopic flux function under local equilibrium. \par 
 In \cite{bmrs4}, we introduce a new approach for the derivation of 
 {\em quenched strong} local equilibrium. 
In the case of subcritical 
hydrodynamic density,  that is,   $\rho(t,x)<\rho_c$,  
we  establish  not only conservation but also  
{\em spontaneous creation} of local equilibrium: 
 we only require starting from a sequence of (possibly deterministic) initial 
configurations with a given macroscopic profile, but with a distribution far away from 
 initial  local equilibrium. 
In the case of supercritical hydrodynamic density $\rho(t,x)>\rho_c$, we prove that 
the local equilibrium property fails, and that, locally around 
``typical points''  of the environment, 
the distribution of the microscopic state is close to the {\em critical measure} 
denoted by $\mu^\alpha_c$: 
this can be viewed as a {\sl dynamic version} of the loss of mass property 
studied in \cite{afgl,bmrs1,bmrs2}. 
 In the case of a critical hydrodynamic density,  that is, 
  $\rho(t,x)=\rho_c$, 
we prove that locally around typical points, the distribution of the system approaches 
the critical measure. This can still be viewed  as a local equilibrium creation result, 
but only in a partial sense. Indeed, outside the situation of an ergodic disorder 
(to which we are not limited), the critical measure may not have critical density, 
nor even any well-defined density.  \par\medskip 
The paper is organized as follows. In Section \ref{sec_cars}, 
we  introduce  our results through an illustrating analysis  of traffic jams. 
In Section \ref{sec_results}, we introduce the model 
and its basic properties. 
Section \ref{sec:convergence} refers to \cite{bmrs1,bmrs2}, 
that is, to the convergence to the critical measure $\mu^\alpha_c$ 
from a supercritical or a critical initial configuration.  
Section \ref{sec:hdl} refers to \cite{bmrs3}, that is, to the 
hydrodynamic limits results (including the supercritical regime). 
Section \ref{sec:loceq} refers to \cite{bmrs4}, that is,
to strong  local equilibrium  in the subcritical and critical regimes 
(creation or conservation), and to 
loss of local equilibrium in the supercritical regime.
\section{A preliminary illustration: traffic jams}
\label{sec_cars}
 In this section we describe some heuristics
for a particular AZRP, where $g(n)=\min(n,1)$ and $p(1)=1$.
This  is a well-known model, namely,
a series of $ M/M/1 $ queues in tandem where each site $x\in\Z$ 
corresponds to a server with service rate $\alpha(x)$.
%
%
Then (see e.g. \cite{par}) provided 
$\lambda<\alpha(x)$ for all $x$, an invariant measure for this process is
the product measure whose marginal at site $x$ is  the geometric distribution with parameter 
$1-\lambda/\alpha(x)$ minus $1$, that has mean value $\lambda/(\alpha(x)-\lambda)$.
The parameter $\lambda$ is the intensity of the Poisson process 
of departures from each queue, 
hence it can be interpreted as the mean current (or \textit{flux}) of customers along the system.
We assume that
\be\label{inf_service}
\forall x\in\Z,\quad 0<c:=\inf_{y\in\Z}\alpha(y)<\alpha(x)
\ee
%
so that the above invariant measure is defined for $\lambda\leq c$.
The value $\lambda=c$ corresponds to what we called ``critical''
 in the introduction.\par\medskip 
This particular AZRP is isomorphic to a 
totally asymmetric simple exclusion process (TASEP),
that models traffic on a one lane highway 
where overtaking is forbidden. In this TASEP, cars (particles) 
are labeled from left to right by integers $n\in\Z$, and
each of them moves one step 
further left  as long as it is not blocked by another car in front 
of it (according to the exclusion rule). 
Particle $n$ jumps with its own 
rate $\alpha(n)$, herafter called its {\em speed}. This 
{\em intrinsic} speed is the one it would reach on an 
otherwise empty road in the absence of any exclusion rule. Whereas 
the original AZRP was endowed with {\em site} disorder, 
the associated TASEP is endowed with {\em particle} disorder. 
Here, \eqref{inf_service} means that we can find cars 
moving at speeds arbitrarily close 
to but not equal to $c$. \par\medskip
The isomorphism between AZRP and TASEP is as follows. 
An AZRP server at site $n\in\Z$ becomes a TASEP car with label $n$,
and customers waiting for this server become vacant TASEP sites between this car
and the next one to the left.
Thus if we denote by $x_n$ 
the position of the $n$-th TASEP particle,
$\eta(n)=x_{n}-x_{n-1}-1$ represents the number of AZRP particles at site $n$. \par
Any time the 
$n$-th TASEP particle jumps to the left (that is $x_n$ 
decreases by $1$), an AZRP customer (i.e. a TASEP hole)  is 
transferred from server $n$ to server $n+1$ at rate 
$\alpha(n)$ (that is $\eta(n)$ decreases by $1$ and $\eta(n+1)$ 
increases by $1$). The displacement of car $n$ corresponds to 
the flux of customers leaving server $n$.
By the isomorphism  between AZRP and TASEP,
invariant  measures can be obtained for the above TASEP: 
under these, inter-particle distances $x_{n+1}-x_n$ are independent 
geometric random variables with  parameter 
$1-\lambda/\alpha(n+1)$, where $0 \leq \lambda \leq c$. \par\medskip 
     \begin{figure}[htbp]
\begin{center}
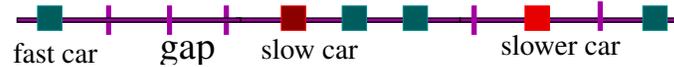
\end{center}
\caption{cars on a highway: fast cars are blue, slow cars are purple, 
and slower cars are red.}
\label{figure:ganaral-cars}
\end{figure}
      \begin{figure}[htbp]
\begin{center}
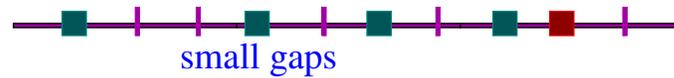
\end{center}
\caption{dense traffic:  all cars are separated by small gaps. }
\label{figure:small-cars}
\end{figure}
          \begin{figure}[htbp]
\begin{center}
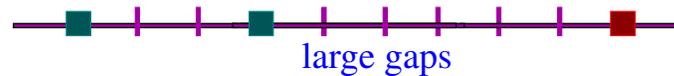
\end{center}
\caption{fluid traffic:  all cars are separated by large gaps, but in front of
 slow cars, the gaps are much larger. }
\label{figure:large-cars}
\end{figure}
\begin{figure}[htbp]
\begin{center}
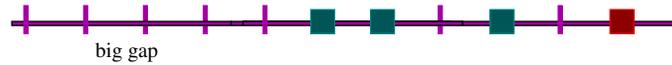
\end{center}
\caption{evolution from fluid traffic}
\label{figure:large-cars2}
\end{figure}
\begin{figure}[htbp]
\begin{center}
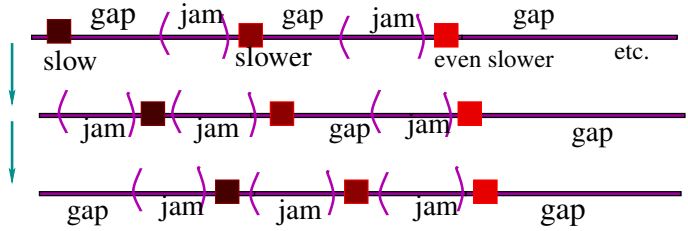
\end{center}
\caption{condensation due to cars regulation: slow cars are dark purple, 
slower cars are light purple and even slower cars are red}
\label{figure:large-cars3}
\end{figure}

 We illustrate these different cases for TASEP  with 
cars with different speeds
(represented by different colors) on a highway, which are now going 
from left to right (see Fig.~\ref{figure:ganaral-cars}). 
When the  traffic is dense,  cars tend to be slowed down by the exclusion rule, 
  compared to which the differences of speeds between them 
  play a lesser role 
  (see Fig.~\ref{figure:small-cars}). 
  There is a phase transition between  dense  and fluid traffic:
  when the  traffic is fluid  (see Fig.~\ref{figure:large-cars}), 
   the difference of speeds plays a dominant role compared to the exclusion rule. 
 There are traffic jams arising behind
 the slowest cars, and big gaps ahead of them and behind the jams 
 generated by the next slower car  (see Fig.~\ref{figure:large-cars2}). 
   Then
  big gaps get reduced  from back to front: the successive jams merge together 
  (see Fig.~\ref{figure:large-cars3}). \par 
   Indeed, to see this, we can define a subsequence  
  $\{\alpha_{n_k}\}$ (of $\{\alpha_n : n \leq 0\}$) of rates 
for successive slower cars
 which strictly decreases to $c$ such that if 
$n_{k+1} < j  < n_k$ then $\alpha_{n_k} \leq \alpha_j $. 
This is done by defining $n_{k+1} = \sup \{ j < n_k: \alpha_j < \alpha_{n_k}\}$ 
and $n_0 =0$.  Let us 
see how cars with labels in the interval $\Z\cap[n_{k+1},n_k]$ 
evolve when we start  with large enough gaps. Neglecting 
the exclusion rule, each car moves at its intrinsic speed 
as long as it is not slowed down by 
the motion of the next car ahead. This generates the following scenario:
\begin{enumerate}   
\item Since car $n_k-1$ has greater speed than car $n_k$, 
its displacement will be bigger, hence the gap between them 
first grows linearly with time, until:
\item in the long run car $n_k$ becomes slowed down by the 
{\em slightly} (if $k$ is large) smaller speed of car $n_{k+1}$ 
and joins the jam already formed behind it by
cars with labels inbetween (which are indeed faster then $n_{k+1}$). 
At this time the existing jam behind car $n_k$ merges with the one 
between $n_{k+1}$;
\item on a longer time scale, this new bigger jam will in turn catch 
up with a slightly slower jam behind car $n_{k+2}$, and so on.
\end{enumerate}   
These steps translate as follows into the AZRP picture, starting 
with large enough occupation numbers (that is, ``supercritical density''):
\begin{enumerate}   
\item The flux from site $n_k -1 $ will be greater than the flux leaving site $n_k$, 
and the occupation number of site $n_k$ will first grow linearly with time, 
which can be interpreted as dynamic condensation; 
\item in the long run the region between sites $n_{k+1}$ and $n_{k}$  
will tend towards an equilibrium with incoming flux close to 
$\alpha_{n_{k+1}}$ on the left end and an outgoing flux 
$\alpha_{n_k}$ close to  matching it, meaning an equilibrium 
at flux $\alpha_{n_{k+1}}$;
\item when the effect 
 from the next block to the left (i.e. from $n_{k+2}$ to $n_k$ )
 reaches $n_k$, the approximate equilibrium in the site 
 interval $\Z\cap[n_{k+2},n_k]$
 will be one for flux $\alpha_{n_{k+2}}$ on the left.
\end{enumerate}   
Since $\alpha_{n_k}\to c$ as $k\to+\infty$,  this suggests 
that as time goes to infinity 
the measure on the AZRP will converge to a 
parametrized by flux $c$.
 A heuristic derivation of the time scale at which 
merging (or \textit{condensation}) occurs, when the service rates 
$\alpha(x)$ are i.i.d. random variables, can be found in \cite{kru}.  \par\medskip
If we zoom out and look at servers from far away (corresponding 
to a scaling limit), one obtains a conservation law with a flux 
function that grows nonlinearly in the ``subcritical'' density range 
(that is, for densities with flux less than $c$), 
and is  truncated at $c$ in the supercritical range.
\par\medskip 
 While we have described the heuristics for the simple 
 model of totally asymmetric case with 
a jump rate one if the site is occupied (for AZRP)
 we obtain the same result for a more complicated 
case of AZRP  (like our general model) where the 
jumps are asymmetric (but not totally asymmetric), 
the jump rate is determined by a function depending on 
the occupation number and the randomness 
of the disorder is much weaker than independent.
\section{Description of the model, basic properties}
\label{sec_results}
In the sequel,  $\mathbb{R}$ denotes the set of real numbers, 
$\mathbb{Z}$  the set of signed integers, $\mathbb{N}=\{0,1,\ldots\}$ 
the set of nonnegative integers
 and $\overline{\mathbb{N}}:=\mathbb{N}\cup\{+\infty\}$. 
For $x\in\mathbb{R}$,  $\lfloor x\rfloor$  
denotes the integer part of $x$, that is largest integer $n\in\mathbb{Z}$ 
such that $n\leq x$.  
The notation $X\sim\mu$ means that a random 
variable $X$ has probability distribution $\mu$.\par 
Let   $\overline{\mathbf{X}}:=\overline{\mathbb{N}}^{\mathbb{Z}}$ 
denote the set of particle configurations, 
and ${\mathbf{X}}:=\mathbb{N}^\mathbb{Z}$ the subset of particle 
configurations with finitely many particles at each site.
 A configuration in $\overline{\mathbf{X}}$ is of the 
 form  $\eta=(\eta(x):\,x\in\mathbb{Z})$ 
 where $\eta(x)\in\overline{\mathbb{N}}$ for each $x\in\mathbb{Z}$. 
 The set  $\overline{\mathbf{X}}$  is equipped with the  coordinatewise  order:  for  
 $\eta,\xi\in\overline{\mathbf{X}}$,  we write 
$\eta\leq\xi$ if and only if $\eta(x)\leq\xi(x)$ for every $x\in\mathbb{Z}$;  in the latter 
inequality, $\leq$ stands for extension to $\overline{\mathbb{N}}$ of
the natural order on $\mathbb{N}$, defined by $n\leq+\infty$ 
for every $n\in\mathbb{N}$, and $+\infty\leq+\infty$. 
  This order is extended to probability measures on 
  $\overline{\mathbf{X}}$: For two 
 probability measures $\mu,\nu$, we write $\mu\le\nu$ if and only if
$\int fd\mu\le\int fd\nu$ for any nondecreasing function $f$ on $\overline{\mathbf{X}}$.  
We denote by $(\tau_x)_{x\in\mathbb{Z}}$ the group of spatial shifts. For $x\in\mathbb{Z}$, 
the action of $\tau_x$  on particle configurations 
is defined by $(\tau_x\eta)(y)=\eta(x+y)$ for every $\eta\in\overline{\mathbf{X}}$, 
 $y\in\mathbb{Z}$. 
Its action on a function $f$ from $\bf X$ or
$\overline\mathbf{X}$ to $\mathbb{R}$ is defined by $\tau_x f:=f\circ\tau_x$.
\subsection{The process and its invariant measures}
\label{subsec:procinv}
Let  $p(.)$ be a probability measure on $\mathbb{Z}$ supported on $\{-1,1\}$.
 We set $p:=p(1)$, $q:=p(-1)=1-p$, and assume $p\in(1/2,1]$,
so that the mean drift of the associated random walk is $p-q>0$.
Let $g:\mathbb{N}\to[0,+\infty)$ be a nondecreasing function such that
\[ 
g(0)=0<g(1)\leq \lim_{n\to+\infty}g(n)=:g_\infty<+\infty\,.
\] 
We extend $g$ to  $\overline{\mathbb{N}}$  by setting $g(+\infty)=g_\infty$.
Without loss of generality, we henceforth assume 
$g(+\infty)=g_\infty=1$. \par
Let $\alpha=(\alpha(x),\,x\in\mathbb{Z})$ (called the environment or disorder) 
be a  $[0,1]$-valued sequence. 
The set of environments is denoted by
 $\mathbf{A}:=[0,1]^{\mathbb{Z}}$.  \par
We consider the 
Markov process  $(\eta_t^\alpha)_{t\geq 0}$
on  $\overline{\mathbf{X}}$  with generator given for any cylinder function  (also called local function, that is depending on finitely many sites) 
$f:\overline{\mathbf{X}}\to\mathbb{R}$  by
\be\label{generator}
L^\alpha f(\eta)  =  \sum_{x,y\in\mathbb{Z}}\alpha(x)
p(y-x)g(\eta(x))\left[
f\left(\eta^{x,y}\right)-f(\eta)
\right]\ee
where, if $\eta(x)>0$,  $\eta^{x,y}:=\eta-\mathfrak{\delta}_x+\mathfrak{\delta}_y$ 
denotes the new configuration obtained from $\eta$ after a particle has 
jumped from $x$ to $y$ (configuration  $\mathfrak{\delta}_x$ has one particle at $x$ and 
no particle elsewhere; addition of configurations is meant coordinatewise). 
 In cases of infinite particle number,  the following interpretations hold:
  $\eta^{x,y}=\eta-\mathfrak{\delta}_x$ if $\eta(x)<\eta(y)=+\infty$ (a particle is removed from $x$), 
  $\eta^{x,y}=\eta+\mathfrak{\delta}_y$ if $\eta(x)=+\infty>\eta(y)$ (a particle is created at $y$),
$\eta^{x,y}=\eta$ if $\eta(x)=\eta(y)=+\infty$. \par\medskip 
For the existence and uniqueness of  $(\eta_t^\alpha)_{t\geq 0}$  
see \cite[Appendix B]{bmrs2}.   
 Recall from \cite{and} that, since $g$ is nondecreasing, 
 $(\eta_t^\alpha)_{t\geq 0}$  
is attractive, i.e.
its semigroup,  denoted by $S(t)$ for $t\geq 0$,  
maps nondecreasing functions (with respect to 
the partial order on $\overline{\mathbf{X}}$) onto nondecreasing functions.
 A graphical construction via a Harris system (\cite{har}) will be a crucial tool
(see e.g. \cite{bmrs2}): attractiveness enables  to construct a completely monotone 
coupling of  a finite number of copies of the process.  \par
The process $(\eta_t^\alpha)_{t\geq 0}$ has the property 
that if $\eta_0\in{\mathbf{X}}$, then almost surely, one has 
$\eta_t\in{\mathbf{X}}$ 
for every $t>0$. In this case, it may be considered as a 
Markov process on $\mathbf{X}$ with generator
(\ref{generator}) restricted to functions $f:{\mathbf{X}}\to\mathbb{R}$. \par
 When the environment $\alpha(.)$ 
is identically equal to $1$, we recover
the {\em homogeneous} zero-range process (see \cite{and} for its detailed analysis). \par\medskip 
For ${\beta}<1$,  
we define the probability measure $\theta_{{\beta}}$ on $\mathbb{N}$
 by \label{properties_a}
\begin{eqnarray}\nonumber 
&&\theta_{{\beta}}(n):=Z({\beta})^{-1}\frac{\beta^n}{g(n)!},\quad n\in\mathbb{N},
\qquad\mbox{where}\quad Z(\beta)
:=\sum_{\ell=0}^{+\infty}\frac{\beta^\ell}{g(\ell)!}\,,\\ \nonumber
&&g(n)!:=\prod_{k=1}^{n} g(k)
\quad\mbox{for}\quad n\in\N\setminus\{0\}, \quad\mbox{and}\quad
g(0)!:=1 \,.
\end{eqnarray}
We denote by  $\mu_\beta^\alpha$  the invariant measure 
of $L^\alpha$ defined (see e.g. \cite{bfl}) as the product measure  
with marginal $\theta_{\beta/\alpha(x)}$ at site $x$:
\be\label{def_mu_lambda_alpha}
\mu^\alpha_{\beta}(d\eta):=\bigotimes_{x\in\Z}\theta_{\beta/\alpha(x)}[d\eta(x)]\,.
\ee
Let
\be\label{inf_alpha}
c:=\inf_{x\in\Z}\alpha(x)\,.
%
\ee
 The measure \eqref{def_mu_lambda_alpha} can be defined on $\overline{\mathbf{X}}$  for
%
${\beta}\in[0,c]$,
%
by using the conventions
\begin{eqnarray}\label{extension_theta}\theta_1&:=&\delta_{+\infty}\,,\\
\label{convention_2}
\frac{\beta}{a}&=& 0\mbox{ if }\beta=0\mbox{ and }a\geq 0\,.
\end{eqnarray}
For $\beta=c$, $\mu^\alpha_c$ is by definition the {\em critical measure}.
The measure \eqref{def_mu_lambda_alpha} is always supported on $\bf X$ if
%
$\beta\in(0,c)\cup\{0\}$.
%
When  $\beta=c>0$, conventions  \eqref{extension_theta}--\eqref{convention_2} 
yield  a measure supported on configurations with infinitely many particles 
at all sites $x\in\Z$ that achieve the infimum in \eqref{inf_alpha}, 
and finitely many particles at other sites.
In particular, this measure is supported on $\bf X$  when the infimum 
in \eqref{inf_alpha} is not achieved. When $c=0$, the measure 
\eqref{def_mu_lambda_alpha} is supported on the empty configuration.\par\medskip 
Since $(\theta_{{\beta}})_{{\beta}\in [0,1)}$ is an exponential family, 
we have a stochastic order relation, that is, 
for $\beta\in[0,\inf_{x\in\mathbb{Z}}\alpha(x)]$, 
 %
$ \mu^\alpha_{{\beta}}$ 
is weakly continuous and stochastically increasing with respect to 
${\beta}\,.$
%
 \subsection{Assumptions on the environment, and consequences}\label{subsec:envt} 
 To state our results,  we introduce two sets of assumptions  on $\alpha$.
The first one says that the environment $\alpha$ has ``averaging properties'', while
the second one sets a restriction
on the sparsity of  ``slow sites'' (where by ``slow sites''  
we mean sites where the disorder variable becomes arbitrarily close 
or equal to the infimum value $c$,  which was  defined by \eqref{inf_alpha}). 
\begin{assumption}\label{assumption_ergo}
 If $c>0$,  
there exists a probability measure $Q_0=Q_0(\alpha)$
on $[0,1]$ such that
\be\label{eq:assumption_ergo}
Q_0(\alpha)=\lim_{n\to+\infty}\frac{1}{n+1}\sum_{x=-n}^0 \delta_{\alpha(x)}
=\lim_{n\to+\infty}\frac{1}{n+1}\sum_{x=0}^n \delta_{\alpha(x)}\,.
\ee
\end{assumption}
\begin{assumption}\label{assumption_dense} The environment $\alpha$ 
has  {\em macroscopically dense defects}, that is, 
there exists a sequence  of sites  $(x_n)_{n\in\mathbb{Z}}$  such that 
\be\label{assumption_afgl}
\forall n\in\mathbb{Z},\,x_n<x_{n+1};\quad
\lim_{n\to\pm\infty}\alpha(x_n)= c \,,
\ee
 and
\[ 
\lim_{n\to\pm\infty}\frac{x_{n+1}}{x_n}=1\,.
\] 
\end{assumption}
Remark that (\ref{assumption_afgl}) implies in particular
\be\label{not_too_sparse}
\liminf_{x\to\pm\infty}\alpha(x)= c\,.
\ee
 By   \eqref{inf_alpha} and \eqref{eq:assumption_ergo},  we have
\be\label{support_Q}
{ c=}\inf_{x\in\mathbb{Z}}\alpha(x)\leq\inf{\rm supp}\,Q_0 \,,
\ee
and  if, for $\beta<1$, we define the mean value of $\theta_\beta$ by
\be\label{mean_theta} R(\beta):=\sum_{n=0}^{+\infty}n\theta_\beta(n)\,,
\ee 
an \textit{average mean density } exists for all $\beta\in[0,c)$: 
\begin{eqnarray}\label{average_afgl}
\overline{R}(\beta)  :=  \lim_{n\to+\infty}
\frac{1}{n+1}\sum_{x=-n}^0 R\left(\frac{\beta}{\alpha(x)}\right)
=\lim_{n\to+\infty}\frac{1}{n+1}\sum_{x=0}^n R\left(\frac{\beta}{\alpha(x)}\right)\,.
\end{eqnarray}
\begin{remark}\label{rk:c>0}
When $c>0$, Assumption \ref{assumption_ergo} is actually {\em equivalent} 
to the fact that the two limits in \eqref{average_afgl} exist and coincide for
every $\beta\in[0,c)$.
\end{remark}
For reasons of concreteness, we denote by  \textbf{(EE)} the stronger  
hypothesis of \textit{environmental ergodicity}: \par\medskip
\noindent
\textbf{(EE)} \textit{The environment has a distribution
 $Q$, for $Q$ a spatially ergodic probability measure on $\bf A$ with marginal $Q_0$,
such that}
$c=\inf{\rm supp}\,Q_0$. \par\medskip 
For instance,   the i.i.d. case  $Q=Q_0^{\otimes\mathbb{Z}}$ satisfies \textbf{(EE)}.  \\
Assumptions \ref{assumption_ergo} and \ref{assumption_dense} are satisfied under 
assumption \textbf{(EE)}; in particular,
$Q$-a.e. realization of the environment $\alpha$  
satisfies Assumption \ref{assumption_dense}. \par\medskip 
By (\ref{average_afgl}), $\overline{R}$ is an increasing $C^\infty$ function 
on $[0,c)$  if $c>0$ (see \cite[Lemma 3.1]{bmrs3}).   
Define the \textit{critical density} by 
\be\label{other_def_critical}
\rho_c:=\overline{R}(c-):=\lim_{\beta\uparrow c}\overline{R}(\beta)\in[0,+\infty]\,,
\ee
and we define $\overline{R}(c):=\overline{R}(c-)=\rho_c$, making 
$\overline{R}(.)$ left continuous at $c$.  If $c=0$, $\overline{R}$ is defined
on $\{0\}$ and equal to $0$ by convention \eqref{convention_2}. In this case, 
we set $\rho_c=0$.\par 
 Moreover, we may define the inverse of $\overline{R}$
on its image $ [0,\rho_c)$,  and, for $\beta\in[0,c)$,  
we reindex the invariant measure $\mu^\alpha_{\beta}$ 
by  setting 
\be\label{invariant_density}
\mu^{\alpha,\rho}:=\mu^\alpha_{\left(\overline{R}\right)^{-1}(\rho)}
\ee   
 The parameter $\rho$ represents the mean particle density in the sense that, for $\mu^{\alpha,\rho}$-almost every configuration $\eta\in\mathbf{X}$, 
\be\label{eq:eta-of-density-rho}
\lim_{n\to+\infty}\frac{1}{n+1}\sum_{x=0}^n\eta(x)=
\lim_{n\to+\infty}\frac{1}{n+1}\sum_{x=-n}^0\eta(x)=\rho\,.
\ee
Note that under \textbf{(EE)},  (\ref{not_too_sparse}) is satisfied,  
$\overline{R}(\beta)$ is well-defined by the ergodic theorem
for all $\beta\in[0,c]$, and we have
\begin{equation}\label{average_afgl_ergodic-and-critical_parameter_max}
\overline{R}({\beta}) 
 =   \int_{(c,1]}R\left[\frac{\beta}{a}\right]dQ_0[a],\qquad
\rho_c  =  \overline{R}(c)\, ,
\end{equation}
\begin{remark}\label{remark_failure}
\begin{itemize}
\item[\textit{(i)}] \ \  Unlike under \textbf{(EE)}, 
 under (\ref{eq:assumption_ergo}), $\lim_{\beta\uparrow c}\overline{R}(\beta)$
could be distinct from $\overline{R}(c)$ as defined by (\ref{average_afgl}) for $\beta=c$.
In fact, the limits in (\ref{average_afgl})
may not exist  for $\beta=c$, or exist and be different
from one another  and/or from $\rho_c$. 
\item[\textit{(ii)}] \ \ In view of \eqref{other_def_critical}, 
it may be tempting to extend the parametrization \eqref{invariant_density} 
to $\rho=\rho_c$
by setting $\mu^{\alpha,\rho_c}=\mu^\alpha_c$. However, this does not 
make real sense, as under \eqref{eq:assumption_ergo}, 
because of \textit{(i)} above, the measure thus 
defined may not satisfy \eqref{eq:eta-of-density-rho} for $\rho=\rho_c$. 
\end{itemize}
\end{remark}
We regard  densities below $\rho_c$ as \textit{subcritical}, density $\rho_c$ as \textit{critical}
and densities above  $\rho_c$ as \textit{supercritical}.
In the former there exists an equilibrium probability measure with this density, 
while we will see that none exists for supercritical densities. 
  \section{Convergence}
\label{sec:convergence}
A natural question for an interacting particle system with multiple equilibria 
is, given a particular equilibrium $\mu$, to determine the set of initial 
configurations for which the particle system converges to
$\mu $. Among classical conservative systems this is fully answered for 
symmetric exclusion processes on $\mathbb{Z}^d $ (see e.g. \cite{ligbook}) 
but in general seems very hard.  
It is often softened to finding large classes of initial distributions 
for which there is convergence.  For asymmetric nearest neighbour exclusion 
processes, the first result of importance was \cite{and2}, see also \cite{mou} 
and \cite{bam} which gave conditions for 
convergence to product measure.  The last is important since it does not require 
an initial translation invariant distribution but only needs a fixed initial 
distribution satisfying the appropriate law of averages.

Theorems \ref{convergence-maximal}  and \ref{prop-nec} below, proved respectively 
in \cite{bmrs1} and \cite{bmrs2},
establish a necessary and sufficient condition for our zero range process to 
converge to $\mu_c^\alpha$.  To our knowledge there are no other comparable 
results for nontrivial assymmetric conservative systems.  Certainly the fact 
(to be justified below), that we consider the maximal equilibrium helps but 
it is to be noted that, \cite{bam} notwithstanding, for the exclusion process 
it is not clear what a necessary and sufficient condition for convergence in 
distribution to $\delta _{\underline 0}$ (and so $\delta _{\underline 1}$) would be
 (where $\underline{0}$ and $\underline{1}$  denote the empty and full configurations,
 defined respectively by
$\underline{0}(x)=1-\underline{1}(x)=0$ for every $x\in\Z$). 
\subsection{Previous results for convergence}
\label{subsec:previous-convergence}
The main paper we refer to is \cite{afgl}, 
which considers  a TAZRP ($p(1)=1$), on $\mathbb{Z}$, 
for  $g(n)=\mathbf{1}_{\{n\ge 1\}}$, with $\alpha(x)\in[c,1]$.
Their results are the following: 
\begin{itemize}
\item[\textit{(a)}] \ \ \cite[Theorem 2.1]{afgl}
For $\alpha$ fixed, the geometric product measures  $\mu_\beta^\alpha$,
for $\beta< \alpha(x)$ for all $x\in\mathbb{Z}$,
 are the extremal invariant measures. The range of the parameter 
 $\beta$ may be either $[0,c)$ (when
$\alpha(x)=c$ for some $x$) or $[0,c]$  (when $\alpha(x)>c$ for all $x$).
\item[\textit{(b)}] \ \ \cite[Proposition 2.2]{afgl}  There does not exist 
any invariant measure with a density above $\rho_c $. 
\item[\textit{(c)}] \ \ \cite[Theorem 2.3]{afgl}
If $\rho_c<+\infty$, if $\eta_0\in\mathbb{N}^\mathbb{Z}$ satisfies 
the supercriticality assumption
\be\label{cond_init-afgl}
\liminf_{n\to\infty}n^{-1}\sum_{x=-n}^0\eta_0(x)>\rho_c\,,
\ee
then $\eta_t^\alpha$ converges in distribution to $\mu^\alpha_c$. 
\end{itemize}
Note that under our assumptions, we obtain $\beta\in[0,c]$ in \textit{(a)}.\par
An immediate consequence of \textit{(c)} is that there is a ``loss of mass'': 
for ``most" of the sites the local density will be $\rho_c$ but our system 
is conservative and the 
hypothesis (\ref{cond_init-afgl}) has a ``global" density (at least 
on the left halfline) strictly higher.

This is an example of {\it condensation} (as exhibited for instance 
in \cite{ev,kru,fk}).  The ``lost" mass accumulates at sites with 
$\alpha$ value close to the minimum $c$
(where ``close" decreases with time). 
If we return to $M/M/1$ queues with service rate $\alpha(x)$ 
at queue $x$, or cars, then,
informally speaking, many clients remain trapped in far 
away slow servers, or lots of space in the road is in front 
of a small number of slow cars.
 \subsection{Our results for convergence}
\label{subsec:convergence}
For our convergence results (from \cite{bmrs1, bmrs2}), 
we also need the following 
 assumption on disorder.  Recall that $c$ is defined by \eqref{inf_alpha}: 
    \begin{assumption}\label{assumption_A1}
    %
We have that $c>0,\, \rho_c<+\infty$,
%
and $\overline{R}$ satisfies the weak convexity assumption\par
{\rm\textbf{(H)}} $\qquad\forall\beta\in[0,c),
\qquad\overline{R}(\beta)-\overline{R}(c)-(\beta-c)\overline{R}^{'+}(c)>0$\par
where $\overline{R}^{'+}(c)$ is the left-hand derivative at $c$ of the convex 
envelope of $\overline{R}$, that is
\[ 
\overline{R}^{'+}(c):=\limsup_{\beta\to c}\frac{\overline{R}(c)-\overline{R}(\beta)}{c-\beta}\,.
\] 
\end{assumption} 
For instance, if $R$ is strictly convex, then for any environment 
satisfying assumption (\ref{average_afgl}), 
$\overline{R}$ is strictly convex
and \textbf{(H)} satisfied. 
A sufficient condition for $R$ to be strictly convex 
(\cite{bs}) is that
$n\mapsto g(n+1)-g(n)$ is a nonincreasing function.
\begin{theorem}
\label{convergence-maximal}
Assume  (\ref{average_afgl}), Assumptions \ref{assumption_dense} 
and  \ref{assumption_A1}.
Then, for all
  $\eta_0\in\mathbb{N}^\mathbb{Z}$ satisfying the supercriticality assumption 
\be\label{cond_init}
\liminf_{n\to\infty}n^{-1}\sum_{x=-n}^0\eta_0(x)\geq \rho_c\,,
\ee
the process  $(\eta_t^\alpha)_{t\geq 0}$ of initial state the configuration 
$\eta_0$ converges in distribution to $\mu_c^\alpha$ when $t\to\infty$.
\end{theorem}
\begin{theorem}\label{prop-nec}
Assume  (\ref{average_afgl}), Assumptions \ref{assumption_dense}
 and  \ref{assumption_A1}.  Assume further that $\eta_0$ satisfies
\[ 
\rho=\liminf_{n\to\infty}n^{-1}\sum_{x=-n}^0\eta_0(x)< \rho_c\,.
\] 
 Then $\eta_t^\alpha$
does not converge in distribution to $\mu_c^\alpha$ as $t\to+\infty$.
\end{theorem}
Together, as claimed, these two results give a necessary and 
sufficient condition for convergence to $\mu_c^\alpha$.\par\medskip
Now note that  Theorem \ref{convergence-maximal} generalises \cite[Theorem 2.3]{afgl}
in the following  ways:\par
\begin{enumerate}
\item the underlying random walk kernel is asymmetric nearest neighbour but 
not necessarily totally asymmetric, 
\item the strict inequality in (\ref{cond_init-afgl}) is replaced by the greater than or equal 
condition of Theorem \ref{convergence-maximal}, 
\item  the special case function $g(n) = \mathbf{1}_{\{n\geq 1\}}$ is removed in favour of any 
$g(.) $ increasing to a finite limit 
 and compatible with Assumption \ref{assumption_A1}.
 \end{enumerate}
In viewing these improvements, one might think that 1) could be improved to at least 
the condition that kernel $p(.)$ is positive mean and of finite range.  In fact, 
surprisingly, this is not possible.  The nearest neighbour requirement is not to 
facilitate the argument.
The reason is that, 
loosely speaking, the result  holds because the incoming flux to the left of the 
origin is the maximal value: $c(p-q)$. This balances the outgoing flux to the right 
of the origin under equilibrium.  How the ``$\liminf $" behaviour for the initial particle 
configuration is achieved is immaterial under nearest neighbour motion and for instance a
configuration that is mostly vacant but which has a sparse set of very high peaks at 
isolated sites poses no problems.    
This is no longer true if the kernel is not nearest neighbour:  an isolated peak with 
a great number of particles initially surrounded by vacant sites may be unable to 
furnish the needed maximal flux.   The following result is established in \cite{bmrs1}.
\begin{theorem}\label{prop_counter}
Assume \eqref{not_too_sparse} and  \eqref{average_afgl}.
Assume further that the jump kernel $p(.)$ is totally asymmetric 
(but not nearest neighbour)  and $p(1)<1$. 
Then there exists $\eta_0\in\N^\Z$ satisfying \eqref{cond_init}, such that 
$\eta^\alpha_t$ does not converge in distribution to $\mu_c^\alpha$ as $t\to +\infty$.
\end{theorem}
The proof of Theorem \ref{convergence-maximal} comes down to showing that any 
limit point of the distributions $S(t)\delta_{\eta_0 }$
must be ``above" and ``below" the target distribution.  The upper bound almost 
follows from \cite{fs} which gives a strong condition for zero range processes 
corresponding to our conditions but with general finite range random walk kernels 
in $\mathbb{Z}^d$ for all $d \geq 1$ so that 
$ \limsup S(t)\delta_{\eta_0 } \ \leq \ \mu_c^\alpha$. 
 Unfortunately  the following  growth condition on $\eta_0$
is imposed in \cite{fs}:
\[
\sum_{n\in\mathbb{N}}e^{-\beta n}\sum_{x:\,|x|=n}\eta_0(x)<+\infty,\quad\forall\beta>0\,.
\] 
 The approach of \cite{bmrs2} is to use (for each $\varepsilon > 0$) a comparison
  with finite Jackson networks on 
intervals $(A_ \varepsilon, a _\varepsilon)$,  where
\begin{eqnarray*}\label{def:ell}
A_\varepsilon:=A_\varepsilon(\alpha)&=&\max\{x\le 0:\alpha(x)\leq c+\varepsilon\}\,,\\
a_\varepsilon:=a_\varepsilon(\alpha)&=&\inf\{x\ge 0:\alpha(x)\leq c+\varepsilon\}\,.\label{def:err}
\end{eqnarray*}
 These networks  evolve according to the  AZRP  rules as if 
the points $A_ \varepsilon$ and $ a _\varepsilon$ were permanently 
occupied by infinitely 
many particles.  As $\varepsilon \rightarrow  0$ the left and 
right endpoints tend respectively 
to $- \infty $ and $+ \infty $.  The nontrivial part was to show 
that this could be done so that: 
\begin{itemize}
\item[\textit{(a)}] \ \
the resulting  ``finite  AZRP"   on 
$(A_ \varepsilon, a _\varepsilon)$ would be 
positive recurrent for each $\varepsilon$; 
\item[\textit{(b)}] \ \ the resulting equilibria converged as $\varepsilon \rightarrow  0$ 
to $\mu_c^\alpha$. 
\end{itemize}
The lower bound is harder to show.  At root it exploits the 
following interface property for one-dimensional nearest neighbour processes 
(which will be substantially deepened in Section \ref{sec:hdl}):
  If two configurations $ \eta_0 $ and  $ \xi_0 $ satisfy the interface condition
\begin{equation}\label{inter-0}
 \exists \ x_0 : \quad \forall \ y \leq x_0 , \ \eta_0(y) \ \leq \ \xi _0(y); \
\forall \ y >x_0 , \ \eta_0(y) \ \geq \ \xi_0(y)
\end{equation}
and if the  AZRP's   $ (\eta_t)_{t\geq 0}$ and  
$ (\xi_t)_{t\geq 0}$ are generated by 
the same Harris system, then the interface property is maintained:
\begin{equation}\label{inter-t}
\forall \ t > 0,\, \exists \ x_t : \quad \forall \ y \leq x_t , \ \eta_t(y) \ \leq \ \xi_t(y); \
\forall \ y >x_t , \ \eta_t(y) \ \geq \ \xi_t(y)\,.
\end{equation}
This simple property permits a recasting of the property of stochastic domination. \par
To show that our process $\eta ^{\alpha}_t $ stochastically dominates $\mu_c^\alpha$
in the limit as $t$ becomes large, it is enough to show that for each 
$\beta<c$, $\eta ^{\alpha}_t $ in the limit dominates  $\mu^{\alpha}_\beta$;
or equivalently (recalling \eqref{invariant_density}) that for each
$\rho<\rho_c$,  
$\eta ^{\alpha}_t $ in the limit dominates  $\mu^{\alpha,\rho}$.  
This will be done if
 (for  $\rho$ fixed)
   $\eta ^{\alpha}_t (x) \geq \xi ^{\alpha, \rho}_t (x)$
   for $x$ fixed
for  $(\xi ^{\alpha, \rho}_s)_{s\ge 0}$  an  AZRP  generated 
with the same Harris system as $(\eta ^{\alpha}_s)_{s\ge 0}$  but initially in 
 $\mu^{\alpha,\rho}$  equilibrium.
We cannot use the interface property to compare $\eta ^{\alpha}_t (x) $ 
and  $ \xi ^{\alpha, \rho}_t (x)$  directly 
so we introduce an intermediary process ($\eta ^{\alpha, t}_s)_{s \geq 0}$ 
that can be compared with both (and indeed everything).
 This process is itself generated by the same Harris system as the 
 other two and is therefore fully defined by specifying
\be\label{source}
 \eta ^{\alpha, t}_0(x)  =(+\infty) {\mathbf 1}_{\{x \leq x_t\}}
\ee
where $x_t $ is a site that is negative, of order $t$,   
with environment value  $\alpha (x_t)$ so that
as $t$ becomes large,  $\alpha (x_t)$ tends to $c$.  
The choice of $x_t$ is not straightforward, it requires assumption \textbf{(H)}. 
 The initial configuration \eqref{source} corresponds to 
placing source/sinks to the left up to site $x_t$ (but since jumps are nearest neighbour,
as seen from the right of $x_t$, this is equivalent to placing 
a source/sink only at $x_t$).  \par\medskip
Due to \eqref{source},  there must be interfaces $I^1 _s $ 
between  $\eta ^{\alpha, t}_s$ and $\eta ^{\alpha}_s$ 
and $I^2 _s $ between  $\eta ^{\alpha, t}_s$ and   $\xi ^{\alpha, \rho}_s$. 
 Our domination result will follow once we have shown 
that with probability tending to one as $t$ becomes large
\begin{itemize}
\item
$I^1 _t $ is highly negative and 
\item $I^2 _t $ is highly positive 
\end{itemize}
as this will imply that around the origin (with high probability)
\be\label{order_interface}
 \eta^{ \alpha }_t(x)  \geq    \eta^{ \alpha, t }_t(x) \ 
 \geq \ \xi ^{\alpha, \rho}_t (x)\,. 
\ee
for all $x$ near the origin.  \par\medskip
To show that e.g. $I^1_t \leq M_t $ (for some well chosen $M_t$) 
it suffices to show that 
on some interval $[N_t,M_t] , \ \sum_{x=N_t} ^ {M_t } (\eta^{ \alpha }_t(x)
 - \eta^{\alpha, t}_t(x)) > 0 $
(and similarly for  $I^2_t \geq M_t $).  We choose $N_t $ and $M_t $ 
to both be of order $t$ 
so that we can use hydrodynamic results to show that at order $t$ 
both $ \sum_{x=N_t} ^ {M_t }  \eta^{ \alpha }_t(x) $ and  
$ \sum_{x=N_t} ^ {M_t }  \eta^{ \alpha, t }_t(x) $ 
are essentially nonrandom and appropriately ordered.   
In fact in this calculation, it was only necessary to understand the hydrodynamic
behaviour and local equilibrium starting from 
 a source initial configuration \eqref{source}. \par 
The key point regarding this behaviour is the following. 
 The hydrodynamic profile created by the source on its right is nonincreasing and
there is critical speed $v_c\geq 0$ such that a front of uniform density $\rho_c$ 
propagates from the source at speed $v_c$. Assumption \textbf{(H)} ensures
that the profile is continuous 
at the end of this front.  This enables us to choose $x_t$ of order 
$-t(v_c+\varepsilon)$ and ensure that \eqref{order_interface} holds on 
$[x_t,0]$ with  any $\rho$ smaller than the hydrodynamic density 
created by the source around the origin, which can be made arbitrarily 
close to $\rho_c$.  \par\medskip 
  A fuller picture of hydrodynamic behaviour 
 and local equilibrium  is discussed in the next sections.
  \section{Hydrodynamics }
\label{sec:hdl}
We begin with the following standard 
definitions in hydrodynamic limit theory. We denote by
$\mathcal M(\mathbb{R})$ the set of Radon measures on 
$\mathbb{R}$. To a particle configuration
 $\eta\in{\mathbf X}$, we associate a sequence of \textit{empirical measures} 
 $(\pi^N(\eta):\,N\in\mathbb{N}\setminus\{0\})$ defined by
\[
\pi^N(\eta):=\frac{1}{N}\sum_{y\in\mathbb{Z}}\eta(y)\delta_{y/N}\in\mathcal M(\mathbb{R})\,.
\]
Let $\rho_0(.)\in L^\infty(\mathbb{R})$, and let 
$(\eta^N_0)_{N\in\mathbb{N}\setminus\{0\}}$ denote 
a sequence of $\bf X$-valued random variables. We say this sequence has 
\textit{limiting density profile} $\rho_0(.)$, if the sequence of empirical measures 
$\pi^N(\eta^N_0)$ converges in probability to the deterministic measure 
$\rho_0(.)dx$ with respect to the topology of vague convergence. \par 
 Let $f:[0,+\infty)\to\R$ be a Lipschitz function, 
and consider the conservation law 
\be\label{conservation_law_0}
\partial_t \rho(t,x)+\partial_x f[\rho(t,x)]=0\,.
\ee
Equation \eqref{conservation_law_0} means that around a point 
where the macroscopic particle density is $\rho$, the instantaneous 
algebraic flux (or current) is 
$f(\rho)$. This equation,
 with given initial condition $\rho_0(.)$ generally does not have 
 strong solutions even if $\rho_0(.)$ is regular, and has infinitely 
 many weak solutions.
However, it has a unique so called {\em entropy} solution, 
that is considered as the {\em physical} solution (\cite{serre}).
 The sequence $(\eta^N_t,\,t\geq 0)_{N\in\mathbb{N}\setminus\{0\}}$ 
 is said to have 
  \textit{hydrodynamic limit} 
$\rho(.,.)$  if: for all 
$t\geq 0,\,(\eta^N_{Nt})_{N\in\mathbb{N}\setminus\{0\}}$ has  
limiting density profile $\rho(t,.)$, which is the 
  \textit{entropy solution}  to the conservation 
  law (\ref{conservation_law_0}) with
 initial datum $\rho_0(.)$. 
\subsection{Previous results for  hydrodynamic limit}
\label{subsec:previous-hydro}
 The hydrodynamic limit of {\em homogeneous} asymmetric zero-range 
process was derived in \cite{av,av1} for step initial conditions
under the assumption of a concave flux function, and in \cite{rez} 
for general initial conditions without the concavity assumption
(the latter result applies to more general attractive models 
with product invariant measures in any space dimension). 
In \cite{lan2}, the hydrodynamic
behavior was studied for an AZRP with a {\em single} 
spatial inhomogeneity exhibiting condensation.\par 
The paper \cite{bfl} derived quenched hydrodynamics 
in the subcritical regime  for an attractive
 AZRP  (i.e. nondecreasing $g(.)$), on $\mathbb{Z}^d$, 
for $p(.)$ finite range, with a disorder such that
$\alpha(x)$ has a finite number of values (thus $\rho_c=+\infty$).\par
In the paper \cite{ks}, hydrodynamics were derived through the 
variational coupling method,  which is effective 
 for a totally asymmetric ZRP ($p(1)=1$), with $g(n)=\mathbf{1}_{\{n\ge 1\}}$, in
 all regimes (subcritical, critical and supercritical).\par\medskip
 The difficulties to prove hydrodynamic limits in our set-up are the following. 
  At supercritical densities, there are no invariant measures;
	moreover we have condensation.
  It is thus impossible to use the traditional approach (\cite{kl}), 
  through block averaging and block estimates,
  since mesoscopic block densities can blow up around condensation sites. 
\subsection{Our results on  hydrodynamic limits}
\label{subsec:hydro}
 The main result of \cite{bmrs3} is the following. 
\begin{theorem}
\label{th_hydro}
(\cite[Theorem 2.1]{bmrs3})
Assume the environment $\alpha$ satisfies  Assumption \ref{assumption_ergo}, 
and the  sequence $(\eta^N_0)_{N\in\mathbb{N}\setminus\{0\}}$ 
has limiting density profile $\rho_0(.)\in L^\infty(\mathbb{R})$. 
For each $N\in\mathbb{N}\setminus\{0\}$, let  
$(\eta^{\alpha,N}_{t})_{t\geq 0}$  denote the process with initial 
configuration $\eta^N_0$ and generator (\ref{generator}). Assume either that 
the initial data is subcritical, that is $\rho_0(.)< \rho_c^{\alpha}$; 
or, that Assumption \ref{assumption_dense} holds.
 Let $\rho(.,.)$ denote the entropy solution to  \eqref{conservation_law_0} 
 with initial datum $\rho_0(.)$,  where $f$ is 
 the flux function defined by  \eqref{def_flux}--\eqref{extension_flux} below. 
Then for any $t>0$,
 the sequence  $(\eta^{\alpha,N}_{Nt})_{N\in\mathbb{N}\setminus\{0\}}$ 
 has limiting density profile $\rho(t,.)$.
\end{theorem}
 To complete the above theorem, we explain how the flux function 
in \eqref{conservation_law_0} is obtained. Let  $x_.:t\mapsto x_t$ 
be  a $\Z$-valued path
representing the position of a moving ``observer'' on the lattice.
We let $\Gamma_{x_.}(t,\eta)$ denote the algebraic current across $x_.$,  
that is the algebraic number of particles crossing the ``observer'' to the right, 
between times $0$ and $t$, when starting from the initial configuration $\eta$. 
In the special case where $x_.$ is identically $0$, 
we obtain the current through the origin between times 0 and $t$, 
simply denoted by $\Gamma_0^\alpha(t,\eta)$: that is, the number of jumps 
from 0 to 1 minus the number of jumps from 1 to 0. 
Then the flux  function $\rho\mapsto f(\rho)$ is defined by
\[
f(\rho):=\lim_{t\to+\infty}t^{-1}\Gamma_0^\alpha(t,\eta^\rho)\,.
\]
We can show that this limit  exists and depends only on $\rho$ 
(and not on the choice of $\eta^\rho$ nor on $\alpha$), where
$\eta^\rho$ denotes a  configuration with density $\rho$, 
that is a configuration  satisfying \eqref{eq:eta-of-density-rho}. 
For $\beta\in[0,c)$, we can compute the stationary current 
under $\mu^\alpha_\beta$, as follows. 
\[
\int_{\mathbf{X}}[p\alpha(x)g(\eta(x))
-q\alpha(x+1)g(\eta(x+1))]d\mu^\alpha_\beta(\eta)=(p-q)\beta\,.
\]
 Recall  the function $\overline{R}$ defined by \eqref{average_afgl}. 
 As a function of the mean  density 
 $\rho=\overline{R}(\beta)$, we  define the flux for our system as
\be\label{def_flux}
f(\rho):=(p-q)\overline{R}^{-1}(\rho)\,.\ee 
And we extend $f$ to $[\rho_c,+\infty)$ by
\be\label{extension_flux}
f(\rho)=(p-q)c,\quad\forall \rho\geq\rho_c\,.\ee
 \begin{remark}\label{remark_flux}
Under Assumption \textbf{(EE)} of an ergodic environment, 
the flux function depends only on the marginal $Q_0$ of the environment.
Under the more general Assumptions \ref{assumption_ergo}--\ref{assumption_dense}, 
the flux function depends on the pair $(Q_0,c)$ 
defined in \eqref{support_Q}--\eqref{average_afgl}. In the latter case, 
the inequality in \eqref{support_Q} can be strict, so $c$ should be regarded
as an additional parameter not contained in $Q_0$, whereas in the ergodic case, 
equality always holds in \eqref{support_Q}.
A simple non-ergodic example is given in Subsection \ref{subsec:ex_hydro} below. 
\end{remark}
To prove Theorem \ref{th_hydro}, we use a reduction principle established 
in \cite{bgrs1,bgrs2,bgrs3,bgrs4,bgrs5} 
for one dimensional conservative attractive processes.
  This method reduces the proof
  of hydrodynamics  for a Cauchy initial condition to that  
  for a Riemann initial condition,  that is, of the form 
\begin{equation}\label{eq:rie}
\rho_0(x)=R_{\lambda,\rho}(x):=\lambda \mathbf 1_{\{x<0\}}+\rho \mathbf 1_{\{x\geq 0\}}
\end{equation} 
for $\lambda,\rho\in\mathbb{R}$. 
  The passage from one to the other
  is similar in spirit to Riemann-based numerical schemes for scalar conservation laws; 
  the difficulty is to control 
  the propagation of the error committed at successive time steps when replacing the actual 
  entropy solution with a piecewise constant
  approximation.  Crucial tools in this reduction are:  
\begin{itemize}
\item[\textit{(i)}]  \ \
The {\em finite propagation} property, that is 
  the fact that discrepancies between two AZRP's, and similarly between two
	entropy solutions of \eqref{conservation_law_0}, propagate with bounded speed. 
\item[\textit{(ii)}] \ \
 The {\em macroscopic stability} property, which states that 
	if two  AZRP configurations  are initially close macroscopically, they remain so at later times.
		In our case this is a consequence of the fact that jumps are 
		nearest neighbour and $g(.)$ is nondecreasing.
\end{itemize}
  To prove hydrodynamics for  the Riemann initial condition  \eqref{eq:rie}, 
  we use a variational characterization  (see e.g. \cite{bgrs5}) of the entropy solution 
 $R_{\lambda,\rho}(t,.)$ of the \textit{Riemann problem} 
 (that is, \eqref{conservation_law_0} with initial datum \eqref{eq:rie}).  Namely, define
	\be\label{variational}
	\mathcal G_{\lambda,\rho}(v)=\left\{
	\ba{lll}
\min_{r\in[\lambda,\rho]}[f(r)-v r] & \mbox{if} & \lambda\leq\rho\,,\\
\max_{r\in[\rho,\lambda]}[f(r)-v r] & \mbox{if} & \lambda\geq\rho\,.
	\ea
	\right.
	\ee
	For such values, $R_{\lambda,\rho}(t,vt)$ is the minimizer (resp. maximizer) 
	in \eqref{variational} if $\lambda\leq\rho$ (resp. $\lambda\geq\rho$). 
	This optimizer is unique for all but countably many values of $v$.\par\medskip
	In order to prove \eqref{variational} at the microscopic level, 
	the main issue is to show that
	\be\label{main_issue}
	\lim_{t\to+\infty}t^{-1}\Gamma_{x_.}(t,\xi^{\alpha,\lambda,\rho})=\mathcal G_{\lambda,\rho}(v)
	\ee
	in probability, where $x_.$ is a path with asymptotic speed $v$, and 
	$\xi^{\alpha,\lambda,\rho}$ is a random configuration
	with profile $R_{\lambda,\rho}(.)$ in \eqref{eq:rie}. Assume for instance $\lambda<\rho$. 
	To define a suitable configuration $\xi^{\alpha,\lambda,\rho}$ in 
	\eqref{main_issue},
	we use a family of AZRP's  $\xi^{\alpha,r,r}_.$, where $r\geq 0$, 
	and $\xi^{\alpha,r,r}_.$ 
	has homogeneous macroscopic density $r$ at time $0$. 
	For $r<\rho_c$, we can choose  $\xi^{\alpha,r,r}_.$  to be 
	an equilibrium process with density $r$ 
	 (that is, with distribution $\mu^{\alpha,r}$, cf. \eqref{invariant_density}). 
	For $r>\rho_c$, such equilibria do not exist. 
Instead, we use what we call ``pseudo-equilibria'', that is, configurations 
$\xi^{\alpha,r,r}$ with a supercritical homogeneous
  macroscopic density profile. We choose $\xi^{\alpha,\lambda,\rho}$ in \eqref{main_issue} 
  as the configuration whose restriction to $x\leq 0$ is $\xi^{\alpha,\lambda,\lambda}$ 
  and whose restriction to $x>0$ is $\xi^{\alpha,\rho,\rho}$. The following main ideas 
  are then involved to derive \eqref{main_issue}:\par\medskip
\begin{itemize}
\item[\textit{(a)}] \ \
 We  prove convergence \eqref{main_issue} for equilibria and pseudo-equilibria, 
	that is $\lambda=\rho=r$.
	This follows from ergodicity in the case $r<\rho_c$ of equilibria, but novel 
	arguments are necessary in the case $r>\rho_c$ of pseudo-equilibria. 
	\item[\textit{(b)}] \ \
 The upper bound in \eqref{main_issue} (proving that the l.h.s. is dominated 
with high probability by the r.h.s.) is the simpler part.
It follows from a coupling argument showing 
that the current in $\xi^{\alpha,\lambda,\rho}_.$ 
cannot exceed  the one in $\xi^{\alpha,r,r}_.$. 
This property can be regarded as ``intuitive''  
because, since $r\in[\lambda,\rho]$, the latter system 
has initially more particles to the left 
and more space to the right. However, mathematically, 
this is related to the macroscopic 
stability property. 
\item[\textit{(c)}]  \ \
 For the lower bound in \eqref{main_issue}, we introduce a novel ``interface process'',
   which gives a more adapted (in our setting) version of the local particle density than 
   the usual block average. This is a random spatially nondecreasing  lattice field 
   $\rho_t(x)$ taking values in $[\lambda,\rho]$ with the following property:  
   in a space region where $\rho_t$ does not fluctuate much, the system is approximately 
   at local equilibrium
  or pseudo-equilibrium, in the sense that $\xi^{\alpha,\lambda,\rho}_t$ is close to 
  $\xi^{\alpha,r,r}_t$ for a random $r=\rho_t(x)$.
	As a nondecreasing function cannot jump too often, for ``most'' 
	values of $v$, using \textit{(a)}, 
	the macroscopic current across a path with velocity $v$ 
	is close to $f[\rho_t(vt)]-v\rho_t(vt)$, 
	which dominates the minimum in \eqref{variational}. 
	Whence the desired lower bound.
\end{itemize}
	The interface process is obtained by looking at simultaneous interfaces of 
	$\xi^{\alpha,\lambda,\rho}_.$ in the sense of (\ref{inter-0})--(\ref{inter-t}) 
	with all equilibria and pseudo-equilibria $\xi^{\alpha,r,r}_.$ for $r\in[\lambda,\rho]$. 
	Precisely, we can define a simultaneous evolution of interfaces $x_.^r$ between 
	$\eta^\alpha_.$ and $\xi^{\alpha,r,r}_.$ so that
	$x_t^r$ is nondecreasing with respect to $r$, and define $x\mapsto\rho_t(x)$ as 
	a generalized inverse of $r\mapsto x_t^r$. 
\subsection{Examples}\label{subsec:ex_hydro}
In this subsection, we illustrate the behaviour of the flux function 
and of solutions to \eqref{conservation_law_0} with some examples.\par\medskip
\runinhead{Dilute limit.} 
It is natural to compare the ``disordered'' flux function 
\eqref{def_flux}--\eqref{extension_flux} to the ``homogeneous'' flux function 
$f_{\rm hom}$ obtained from the same AZRP in a homogeneous 
environment $\alpha_{\rm hom}(x)\equiv 1$. 
In general, there is no simple relation between these two fluxes. 
However, the relation becomes simple and natural in the so-called 
{\em dilute limit} where $\alpha(x)=1$ at ``most'' sites.
More precisely, let $Q_0$ be the probability measure in  Assumption \ref{assumption_ergo}, 
\textbf{(EE)} and \eqref{average_afgl_ergodic-and-critical_parameter_max}.
We consider an i.i.d. random environment in which, for each $x\in\Z$, $\alpha(x)$ 
is chosen chosen according to $Q_0$ with probability $\varepsilon\geq 0$. and 
equal to $1$ with probability $(1-\varepsilon)$. Thus $\alpha(x)$ has distribution
\be\label{def_qeps}
Q_0^\varepsilon:=(1-\varepsilon)\delta_1+\varepsilon Q_0\,.
\ee
The value $\varepsilon=0$ corresponds to the homogeneous environment 
$\alpha_{\rm hom}$, while $\varepsilon=1$ corresponds to the i.i.d.
environment with marginal $Q_0$. The dilute limit is the limit 
$\varepsilon\to 0+$.
Let $f^\varepsilon$ denote the flux function \eqref{def_flux}--\eqref{extension_flux} 
produced by the environment with marginal \eqref{def_qeps}.
For $\varepsilon=0$, $f^0=f_{\rm hom}$ is the flux function 
for the homogenous AZRP. 
It follows from \eqref{def_flux}--\eqref{extension_flux} that
\be\label{flux_hom}
f_{\rm hom}(\rho)=(p-q)R^{-1}(\rho)\,.
\ee
For $\varepsilon\in(0,1]$, there is no simple relation 
between $f^\varepsilon$ and $f_{\rm hom}$.
However, using \eqref{def_flux}--\eqref{extension_flux} 
and \eqref{average_afgl_ergodic-and-critical_parameter_max}, we can show that
\begin{eqnarray*}\label{dilute_critical_4}
\lim_{\varepsilon\to 0}f^\varepsilon(\rho) & = & f_d(\rho)
\end{eqnarray*}
where $f_d$ is the dilute limit of the flux function, defined by
\be\label{dilute_flux}
f_d(\rho):=\left\{
\ba{lll}
f_{\rm hom}(\rho)
& \mbox{if} & \rho<\rho_c\\ 
 (p-q) c & \mbox{if} & \rho\geq\rho_c
\ea
\right\}=f_{\rm hom}(\rho)\wedge  (p-q) c\,.
\ee
This limit can be understood 
intuitively as follows. As $\varepsilon\to 0$,  slow sites 
are very rare, hence the system exhibits long homogeneous stretches 
where it behaves as a homogeneous process. Thus the memory 
of slow sites is only retained by the flux truncation, but not 
by the shape of the flux function prior to truncation. Note in 
particular that only the infimum of the support of $Q_0$
in \eqref{def_qeps} (but not details of the distribution) 
is involved in the dilute limit.  \par\medskip
We defined above the dilute limit from the limit of a sequence of random environments. 
Another point of view is to construct a {\em single deterministic} environment 
equal to $1$ except on a $0$ density subset of $\Z$, where it asymptotically 
approaches its infimum value $c$. Precisely, let $\alpha(.)$ be an environment 
satisfying the conditions of Assumption \ref{assumption_dense},  with values in $[c,1]$, 
and such that 
\[
\alpha(x)=
\ba{lll}
1 & \mbox{if} & x\not\in\{x_n:\,n\in\Z\}\,.
\ea
\]
 Assume moreover $\lim_{n\to\pm\infty}(n/x_n) =  0$. 
For such an environment, the flux function defined by 
\eqref{def_flux}--\eqref{extension_flux} is exactly $f_d$ 
given by \eqref{dilute_flux}. As announced in Remark \ref{remark_flux}, 
this example shows that, outside the case of an ergodic random environment, 
the flux function is not entirely determined by the empirical
distribution $Q_0$ in \eqref{assumption_afgl}. Indeed, in this case we 
have $Q_0=\delta_1$, which does not give any information on $c$. \par\medskip
\runinhead{Supercritical entropy solutions.}
We now describe the consequences of the flat line \eqref{extension_flux} 
on the behaviour of entropy solutions through the analysis of the 
so-called {\em Riemann} problem with
initial data of the form  \eqref{eq:rie}.  The following result can be 
obtained (see \cite{bmrs3}) using \eqref{variational}. 
\begin{proposition}\label{cor_riemann}
 Assume $+\infty>\lambda\geq\rho_c>\rho$. 
 Let
\be\label{critical_speed}
 v_c(\rho):=\inf_{r\in[\rho,\rho_c)}\frac{f(\rho_c)-f(r)}{\rho_c-r}  =
\inf_{r\in[\rho,\rho_c)}\frac{\widehat{f}(\rho_c)-\widehat{f}(r)}{\rho_c-r}
=\widehat{f}\,'(\rho_c-) 
\ee
 where $\widehat{f}$ denotes the concave envelope of $f$ on $[\rho,\rho_c]$. 
 In particular, if $f$ is concave,
\[ 
v_c(\rho)=f'(\rho_c-)=\left\{\int_{[c,1]}\frac{1}{a}R'\left[
\frac{c}{a}
\right]
dQ_0(a)\right\}^{-1}\,.
\] 
Then, for every $t>0$, we have
\begin{eqnarray}
\label{riemann_left}
R_{\lambda,\rho}(x,t) & = & \lambda,\quad\forall x<0\\ 
R_{\lambda,\rho}(x,t) & = & R_{\rho_c,\rho}( x,t ),\quad\forall x>0 \label{riemann_right}\\ 
 \lim_{t\to+\infty}R_{\lambda,\rho}(x,t) & = & \rho_c,\quad\forall x\geq 0\label{endupcrit}\\
R_{\lambda,\rho}(0+,t) & = & \rho_c\label{critical_origin}\\ 
R_{\lambda,\rho}(x,t) & = & \rho_c,\quad\forall x\in(0,tv_c(\rho))\label{front}\\
R_{\lambda,\rho}(x,t) & < & \rho_c,\quad\forall x>tv_c(\rho)\label{nofront}
\end{eqnarray}
\end{proposition}
Property \eqref{riemann_left} states that
the initial  constant  density 
is not modified to the left of the origin.
This is related to the fact that $f$ 
is nondecreasing, hence characteristic velocities are always nonnegative. \par
Property \eqref{riemann_right} states that the solution to the right of 
the origin does not depend on the supercritical initial density $\lambda$ 
on the left. Formally,  we may thus consider that this is also the solution 
for ``$\lambda=+\infty$'', which corresponds to placing sources to the left 
of the origin. In particular, for $\rho=0$, we recover the source solution 
used in Subsection \ref{subsec:convergence}, that is the hydrodynamic limit
starting from the particular source configuration \eqref{source}.\par 
Properties \eqref{endupcrit}, \eqref{critical_origin}, \eqref{front} 
are signatures of the phase transition.
They express the fact that, regardless of the supercritical value on the left side,
supercritical densities are blocked, and
the right side is dominated by the critical density.
In particular, \eqref{front}--\eqref{nofront} state that
a front of critical 
density propagates to the right from the origin  at speed $v_c$ if $v_c(\rho)>0$. \par\medskip
\runinhead{An explicit case:  $M/M/1$ queues in series.}
Consider the $M/M/1$ queues in series, 
that is $g(n)=\min(n,1)$, with total asymmetry ($p=1$, $q=0$),  cf. Sect. \ref{sec_cars}. 
With this choice of $g$,
\eqref{mean_theta} and \eqref{flux_hom} write
\be\label{mm1}
R(\beta)  =  \frac{\beta}{1-\beta},\quad
f_{\rm hom}(\rho) = \frac{\rho}{1+\rho}\,.\label{mm1_fhom}
\ee
In view of \eqref{mm1_fhom}, the dilute limit $f_d$ defined 
in \eqref{dilute_flux} writes here
\be\label{mm1_dilute}
 f_d(\rho)=\left[\frac{\rho}{1+\rho}\right]\wedge c=\left\{
\ba{lll}
\displaystyle{\frac{\rho}{1+\rho} } & \mbox{if} & \displaystyle{\rho<\rho_c:=\frac{c}{1-c}}\,,\\ 
c & \mbox{if} & \rho\geq\rho_c\,.
\ea
\right.
\ee 
 Since $f_d$ defined by \eqref{mm1_dilute} is concave, \eqref{critical_speed} yields 
\[ 
v_c(\rho) =f'_{\rm hom}(\rho_c^-)=(1-c)^2\,.
\] 
 \section{Local equilibrium}\label{sec:loceq}
 We  now come to results on 
strong local equilibrium with respect to the hydrodynamic 
limit  (\ref{conservation_law_0}).\par 
A natural extension of the convergence theorems is to establish results for the following 
question. Let us fix a realization of the environment $\alpha (.) $ and 
suppose given a sequence of initial configurations  $\{\eta^{N}_0\}_{N\ge 0}$ 
which correspond 
to a ``profile" $ \rho_0 $ in the sense that, for every $a,b\in\R$,
\begin{equation}  \label{HIN}
 \frac{1}{N} \sum_{aN < k < bN} \eta^{N}_0 (k)  \ - \ \int _a ^b \rho_0(x) dx \ \rightarrow \ 0\,.
\end{equation}
Let $x_N \ \in \ \mathbb{Z} $ satisfy  $N^{-1}x_N \ 
\rightarrow \ u \ \in \  (- \infty, \infty ) $. 
What can be said about the behaviour of  $\eta^{\alpha,N}_{Nt} $  around $x_N $?\par
We suppose that the entropy solution  $ \rho(.,.) $  to the associated hydrodynamic 
equation (\ref{conservation_law_0}) is continuous at  $(t,u)$. 
 It is reasonable to believe that  $ \tau_{[x_N]}\eta^{\alpha,N}_{Nt}
$  ``looks" like the equilibrium corresponding to  $ \rho(t,u) $ 
(with environment $\alpha$ suitably shifted).
This question is anticipated by several works on conservative systems without disorder. 
In translation invariant cases where there is a family of equilibria 
 $\{\mu^r\}_{r\geq 0} $  
and the initial configurations are random
there is \textit{{conservation of local equilibrium}} if, when 
\[
\lim_{N\to\infty} \tau_{[uN]}\eta^N_0 = \mu^{\rho_0(u)}
\]
in distribution for each  continuity point  $u$ of $\rho_0(.)$ 
and for each $(t,u)$  as above,  then  
\begin{equation}  \label{HLIM}
\lim_{N\to\infty} \tau_{[uN]}\eta^N_{Nt} =  \mu^{\rho(t,u)}
\end{equation} 
in distribution.  In \cite{lan}, conservation of local 
equilibrium was proved for a 
homogeneous zero-range process with a strictly convex flux, 
under an initial invariant product measure.
Following \cite{and2}, 
 \cite{bam} showed for finite range nonzero mean exclusion processes, that (\ref{HLIM})
held with no further assumptions on $\eta^N_0$ beyond the profile hypothesis (\ref{HIN}). 
 We call this  a ``spontaneous creation of local equilibrium.''  \par\medskip
%
%
In our family of models it is natural to hope for similar results.  
We discuss three differences:
\begin{enumerate}
\item 
Our system is not translation invariant and in fact  we cannot expect 
a limit as in (\ref{HLIM}) because the environment varies as $N\to+\infty$:
  given appropriate  $(t,u)$ 
it is possible to find sequences $x_N $ and $y_N$ that converge
 macroscopically to $x$ so that 
$ \Exp( f( \tau_{x_N}  \eta^{\alpha,N}_{Nt})) \ - \ \Exp( f( \tau_{y_N}  \eta^{\alpha,N}_{Nt}))$ 
 does not converge to zero as $N$ becomes large. 
\item  
Around a fixed point, say 0, the environment is fixed and will satisfy 
$\alpha (x) < c $ for  $x$  close to the point, however when considering 
$x_N $ which (in scale $N$)  converges  to, say the origin, 
it is perfectly possible that  $\alpha ( x_N) $
converges to $c$.
\item 
The value   $\rho(t,u) $  may be supercritical, that is strictly greater than
 $\rho_c $ the maximum density for an equilibrium.
 \end{enumerate}
The first point is dealt with by a simple reformulation of the result but 
the second and third questions require a more substantive response.  
The second requires us to consider  AZRP  for which the 
occupation number  $+\infty $ is permitted.  
For the third, as with other papers analyzing condensation phenomena, 
we expect that the limiting density will be $\rho_c $ and not  $\rho(t,u)$. \par\medskip
Our first results concern cases where the entropy solution  $\rho(.,.) $ 
is continuous at  $(t,u)$  and has a value strictly below $\rho_c$. 
\begin{theorem}\label{th_strong_loc_eq}
Under assumptions of Theorem \ref{th_hydro}, 
the following holds for every $(t,u)\in(0,+\infty)\times\mathbb{R}$:
let $\psi:{\mathbf X}\to\mathbb{R}$  
be a bounded local function, and $(x_N)_{N\in\mathbb{N}}$ a sequence  of sites 
such that  $u=\lim_{N\to+\infty}N^{-1}x_N$. Then 
 if $\rho(.,.)$ is continuous at $(t,u)$ and $\rho(t,u)<\rho_c$, 
\be\label{eq:loc_eq}
\lim_{N\to+\infty}\left[
\Exp 
\psi\left(
\tau_{x_N}\eta^{\alpha,N}_{Nt}
\right)-\int_{\mathbf X}\psi (\eta)d\mu^{
\tau_{x_N}\alpha,\rho(t,u)
}(\eta)
\right]=0\,.
\ee
\end{theorem}
This result relied upon coupling.  It is sufficient to show the desired convergence 
for increasing cylinder functions $\psi$.  So it is sufficient to show that for 
 $  r \ < \ \rho(t,u) $ we have 
\[
\liminf_{N \rightarrow \infty }
\Exp 
\psi\left(
\tau_{x_N}\eta^{\alpha,N}_{Nt}
\right)-\int_{\mathbf X}\psi (\eta)d\mu^{
\tau_{x_N}\alpha,r
}(\eta)
 \ \geq \ 0
\] 
for $N$ large; and similarly for  $r \ > \ \rho(t,u) $.   
So it will be sufficient to show that for $N$ large in a 
$N$ order interval centered around $x_N$, the process  $\eta^{\alpha,N}_{tN}$  
dominates  $\xi^{\alpha, r}_{Nt} $ where as before  $\xi^{\alpha, r}_{.t} $
is a  AZRP  run with the given Harris system in $\mu^{\alpha, r }$ equilibrium.\par 
Our approach shadows that of \cite{bam} but is a bit more complicated than 
this argument which dealt with the exclusion process.  In particular 
the key tool of \cite{kos},  which relies on global strict concavity 
or convexity of the flux function,  is no longer available here.  \par
Under the general assumptions on function $g(.)$, the flux function $f$ need 
not be globally convex or globally concave.  Nonetheless away from the critical 
value the function $f$ is analytic and therefore in any interval of densities, 
there will exist subintervals on which the flux is either globally convex or 
globally convex.  Thus given  $ r < \rho(t,u) $ we can find $ r < r_1 < r_2 < \rho (t,u) $ 
so that $f$ is either convex or concave on  $(r_1, r_2)$.   
For concreteness suppose it is 
concave, then we may ``impose" a system of \textit{priorities} 
(or \textit{classes}, cf. \cite{ak}) on 
 $\eta^{\alpha,N}_{s} $  particles that are not coalesced 
with  $\xi^{\alpha, r}_{s} $  particles.  
This means that higher priority particles will be faster than lower priority ones. 
So a given small, but order $N$, interval $I$ can be placed between two similar intervals 
$J_l $ and $J_r $.
In a short time (unless they become coalesced) fast particles in $J_l $ will overtake 
slow particles in $J_r $.  
This entails coalescence of many uncoalesced   $\xi^{\alpha, r}_{s} $ 
particles originally in $I$.  
This argument is repeated until in an interval around $x_N $, 
the ``density" of  uncoalesced  
 $\xi^{\alpha, r}_{s} $  particles will be very small. 
 We  can then argue as in \cite{bam} to show the required domination. \par\medskip
Next, we consider points at which the hydrodynamic density is supercritical 
 or critical. As discussed above,  there does not exist a corresponding 
 equilibrium measure,
so one cannot expect the same type of convergence as above. 
We expect and to some extent show 
one has local convergence to the {\em critical} 
quenched invariant  measure.  
 We now introduce an additional assumption that is needed. 
\begin{definition}\label{def_typical}
Let $(x_N)_{N\in\mathbb{N}}$ be a sequence of sites such that 
$N^{-1}x_N$ converges to $u\in\mathbb{R}$ as $N\to+\infty$.  
The sequence $(x_N)_{N\in\mathbb{N}}$ is {\em typical} if and only if 
any subsequential limit $\overline{\alpha}$ of 
 $(\tau_{x_N}\alpha)_{N\in\mathbb{N}}$ 
has the following properties:
\begin{itemize}
\item[\textit{(i)}]\ \ For every $z\in\mathbb{Z}$, 
$\overline{\alpha}\in{\mathbf{B}}:=(c,1]^\mathbb{Z}$.
\item[\textit{(ii)}]\ \ 
$\liminf_{z\to-\infty}\overline{\alpha}(z)=c$.
\end{itemize}
\end{definition}
 The interpretation of the word ``typical'' is the following. When $\alpha$ does 
 not achieve its infimum value $c$ (that is $\alpha\in\mathbf{B}$),
\textit{(i)--(ii)} says that the environment as seen from $x_N$ 
shares  key properties of the environment 
seen from a fixed site (say the origin). 

The typicality assumption will in fact be needed only to show that the  
microscopic distribution is locally dominated by the critical measure.
We observe however that even in the critical case, where the invariant measure does exist, 
we are not able to prove the statement without the typicality assumption. 
\begin{theorem}\label{th_strong_loc_eq_super}
Under assumptions  of Theorem \ref{th_hydro}, 
the following holds for every $(t,u)\in(0,+\infty)\times\mathbb{R}$.
Let $\psi:{\mathbf X}\to\mathbb{R}$ 
be a bounded local function, and $(x_N)_{N\in\mathbb{N}}$ a typical sequence  of sites 
such that  $u=\lim_{N\to+\infty}N^{-1}x_N$. Then: 
\begin{itemize}
\item[(i)]\ \ If $\psi$ is nondecreasing, 
\be\label{eq:loc_eq_super}
\lim_{N\to+\infty}\left[
\Exp 
\psi\left(
\tau_{x_N}\eta^{\alpha,N}_{Nt}
\right)-
\int_{\overline{\mathbf X}}\psi(\eta)d\mu_c^{
\tau_{x_N}\alpha
}(\eta)
\right]^+=0\,.
\ee
\item[(ii)]\ \ If  $\rho_*(t,u)\geq\rho_c$ 
(where $\rho_*(t,u):=\liminf_{(t',u')\to(t,u)}\rho(t',u')$), then
\be\label{eq:loc_eq_super_bis}
\lim_{N\to+\infty}\left[
\Exp 
\psi\left(
\tau_{x_N}\eta^{\alpha,N}_{Nt}
\right)-
\int_{\overline{\mathbf X}}\psi(\eta)d\mu_c^{
\tau_{x_N}\alpha
}(\eta)
\right]=0\,.
\ee
\end{itemize}
\end{theorem}
The argument we give for the upper bound is now close to that for the 
convergence upper bound: we show that we can find $y_N < x_N < z_N $ so that
the ``finite"  AZRP  on $(y_N, z_N)$ with infinitely many 
particles at $y_N$ and $z_N $ 
has an equilibrium very close to  $\mu^{\alpha}_c$.  The lower bound is 
essentially supplied by the arguments for Theorem \ref{th_strong_loc_eq}.\par\medskip
 Without the above typicality 
assumption, we are not able to prove such a statement at given times, but we can obtain 
a weaker time-integrated result. 
\begin{theorem}\label{prop_strong_loc_eq_super_cesaro}
 Under assumptions and notations of Theorem \ref{th_hydro}, the following holds.
Let $t>0$, and  $(x_N)_{N\in\mathbb{N}}$ be an {\em arbitrary} sequence of sites such that 
$\lim_{N\to+\infty}N^{-1}x_N=u$, where $u\in\mathbb{R}$ 
 is such that $\rho_*(t,u)\geq\rho_c$.   Let $\psi:\overline\mathbf{X}\to\mathbb{R}$ be a 
 continuous local function. Then
\be\label{eq:loc_eq_super_cesaro}
\lim_{\delta\to 0}\limsup_{N\to+\infty}\left|
\frac{1}{\delta}\int_{t-\delta}^{t}\Exp 
\psi\left(
\tau_{x_N}\eta^{\alpha,N}_{Ns}
\right)ds-\int_{\overline{\mathbf X}}\psi(\eta)d\mu^{
\tau_{x_N}\alpha
}_c(\eta)
\right|=0\,.
\ee
\end{theorem}
\begin{remark}\label{remark_whereas}
Whereas in Theorems \ref{th_strong_loc_eq} and \ref{th_strong_loc_eq_super}, 
the test function $\psi$ is defined on $\mathbf{X}$, in Theorem 
\ref{prop_strong_loc_eq_super_cesaro}, it is defined on $\overline{\mathbf{X}}$. 
Indeed in Theorem \ref{th_strong_loc_eq}, the fact that $\rho(t,u)<\rho_c$
implies that the subcritical measure $\mu^{\tau_{x_N}\alpha,\rho(t,u)}$ in 
\eqref{eq:loc_eq} is supported on $\mathbf{X}$. In Theorem \ref{th_strong_loc_eq_super},
statement \textit{(i)} in Definition \ref{def_typical} of the typicality assumption ensures 
that the critical measure $\mu_c^{\tau_{x_N}\alpha}$ 
in \eqref{eq:loc_eq_super}--\eqref{eq:loc_eq_super_bis} is still supported on $\mathbf{X}$. 
However in Theorem  \ref{prop_strong_loc_eq_super_cesaro}, owing to the absence of 
the typicality assumption, the critical measure $\mu_c^{\tau_{x_N}\alpha}$ 
in \eqref{eq:loc_eq_super_cesaro} may have infinitely many particles at sites 
$x\in\Z$ such that $\alpha(x)=c$.
\end{remark}
Our approach is to first note that, again, the lower bound is achieved 
via Theorem \ref{th_strong_loc_eq}. 
The upper bound is more intricate.  We argue that in the supercritical 
regime the current through any point must (``on average")  equal $c$. 
We then argue that any 
invariant measure that dominates  $\mu^{\alpha}_c$  and has current $c$ 
must equal  $\mu^{\alpha}_c$.   The (unfortunately necessary) Ces\`aro means 
are to give an invariant measure for any limit measure and to ensure that 
the limit measure must have flux $c$. \par\medskip
\runinhead{Example.} 
An illuminating particular case of the above theorems is when the initial datum is uniform.
Assume 
\be\label{initial_uniform}\rho_0(x)\equiv\rho\ee
for some $\rho\geq 0$. Then $\rho(t,x)\equiv\rho$ for all $t>0$.
Specializing Theorems \ref{th_strong_loc_eq} and \ref{th_strong_loc_eq_super} 
to $u=0$, we obtain that $\eta_{Nt}^{\alpha,N}$ converges in distribution to 
$\mu^{\alpha,\rho}$ if $\rho<\rho_c$, or to $\mu_c^\alpha$ if $\rho\geq\rho_c$. 
We may in particular achieve (\ref{initial_uniform}) as follows by a sequence 
of initial configurations $\eta^N_0=\eta_0$ independent of $N$.
\begin{itemize}
\item[\textit{(i)}] \ \ \textit{Stationary initial state}. 
Let $\eta^N_0=\eta_0\sim\mu^{\alpha,\rho}$, 
with $\rho<\rho_c$. Since $\mu^{\alpha,\rho}$ is an invariant measure 
for the process with generator (\ref{generator}), for every $t>0$, we have 
$\eta^{\alpha,N}_{Nt}\sim\mu^{\alpha,\rho}$.  As a result, the expression between 
brackets in (\ref{eq:loc_eq}) vanishes for every $N\in\mathbb{N}\setminus\{0\}$ and 
$t\geq 0$. Hence there is {\em conservation} of local equilibrium (since 
(\ref{eq:loc_eq}) already holds for $t=0$), but there is in fact nothing to prove, 
since this conservation follows form stationarity.\par \smallskip
\item[\textit{(ii)}] \ \ \textit{Deterministic initial state}. Let 
$\eta^N_0=\eta_0$, 
where $\eta_0\in{\mathbf{X}}$ has density $\rho$,  cf. \eqref{eq:eta-of-density-rho}. 
This 
implies that the sequence $(\eta^N_0)_{N\in\mathbb{N}\setminus\{0\}}$ has density 
profile $\rho_0(x)\equiv\rho$.
 Here Theorems \ref{th_strong_loc_eq} and \ref{th_strong_loc_eq_super} 
are no longer void statements as in \textbf{(i)}. When specialized to $u=0$, 
they yield a large-time convergence result for our process. Namely,
Theorem \ref{th_strong_loc_eq} implies that $\eta^\alpha_t\to\mu^{\alpha,\rho}$ 
in distribution if $\rho<\rho_c$, and Theorem \ref{th_strong_loc_eq_super} 
implies $\eta^\alpha_t\to\mu^{\alpha}_c$ if $\rho\geq\rho_c$. 
\end{itemize}
Pushing the analysis  of \textit{(ii)}  further, 
  we can  show that the first limit in 
  \eqref{eq:eta-of-density-rho} is irrelevant, and  
  derive from Theorems \ref{th_strong_loc_eq} and 
  \ref{th_strong_loc_eq_super} the following convergence result.
\begin{theorem}\label{cor_dream_bfl}
Let $\eta_0\in{\mathbf X}$ be such that,  for some $\rho>0$, 
\[ 
\lim_{n\to+\infty}\frac{1}{n}\sum_{x=-n}^0\eta_0(x)=\rho\,.
\] 
Then $\eta_t^\alpha$ converges in distribution as $t\to+\infty$ 
to $\mu^{\alpha,\rho\wedge\rho_c}$.
\end{theorem}
\begin{remark}\label{remark_improvement} 
\begin{itemize}
\item[\textit{(i)}]\ \ In the case $\rho<\rho_c$, 
 Theorem  \ref{cor_dream_bfl} solves the  convergence  
problem posed in  \cite[pp 195-196]{bfl}.\par\noindent  
\item[\textit{(ii)}]\ \ In the case $\rho\geq\rho_c$, 
 Theorem  \ref{cor_dream_bfl} is a partial improvement over 
 Theorem \ref{convergence-maximal}.
It improves the latter in the sense that we do not need a 
weak convexity assumption on $g$, 
but is it less general with respect to initial conditions.
\end{itemize}
\end{remark}
\begin{acknowledgement}
C.B. and E.S. thank respectively the organizers of PSPDE V and PSPDE VII 
for their invitation to  participate  in  these conferences, 
and therefore giving them the opportunity to present their work, 
then to write the present text. \par
 This work was partially supported by laboratoire MAP5,
 grants ANR-15-CE40-0020-02 and ANR-14-CE25-0011,
 Simons Foundation Collaboration grant 281207 awarded to K. Ravishankar.
 This work has been conducted within the FP2M federation (CNRS FR 2036).
 C.B., T.M. and K.R. thank
 Universit\'{e} Paris Descartes for hospitality.
\end{acknowledgement}
%
%
%
%

%
\end{document}